\documentclass{article}

\usepackage{graphicx}
\usepackage{color}
\usepackage{indentfirst}
\usepackage{amsmath,amsfonts,amsthm,amssymb,amsrefs}
\usepackage{mathrsfs}
\usepackage{amscd}
\usepackage{hyperref}

\DeclareMathAlphabet{\mathsfsl}{OT1}{cmss}{m}{sl}

\newtheorem{thm}{Theorem}[section]

\newtheorem{cor}[thm]{Corollary}
\newtheorem{prop}[thm]{Proposition}

\newtheorem*{thm*}{Theorem}

\theoremstyle{definition}

\newtheorem{rem}[thm]{Remark}

\begin{document}

\def\G{{\Gamma}}
  \def\d{{\delta}}
  \def\ci{{\circ}}
  \def\e{{\epsilon}}
  \def\l{{\lambda}}
  \def\L{{\Lambda}}
  \def\m{{\mu}}
  \def\n{{\nu}}
  \def\o{{\omega}}
  \def\O{{\Omega}}
  \def\Th{{\Theta}}\def\s{{\sigma}}
  \def\v{{\varphi}}
  \def\a{{\alpha}}
  \def\b{{\beta}}
  \def\p{{\partial}}
  \def\r{{\rho}}
  \def\ra{{\rightarrow}}
  \def\lra{{\longrightarrow}}
  \def\g{{\gamma}}
  \def\D{{\Delta}}
  \def\La{{\Leftarrow}}
  \def\Ra{{\Rightarrow}}
  \def\x{{\xi}}
  \def\c{{\mathbb C}}
  \def\z{{\mathbb Z}}
  \def\2{{\mathbb Z_2}}
  \def\q{{\mathbb Q}}
  \def\t{{\tau}}
  \def\u{{\upsilon}}
  \def\th{{\theta}}
  \def\la{{\leftarrow}}
  \def\lla{{\longleftarrow}}
  \def\da{{\downarrow}}
  \def\ua{{\uparrow}}
  \def\nwa{{\nwtarrow}}
  \def\swa{{\swarrow}}
  \def\nea{{\netarrow}}
  \def\sea{{\searrow}}
  \def\hla{{\hookleftarrow}}
  \def\hra{{\hookrightarrow}}
  \def\sl{{SL(2,\mathbb C)}}
  \def\ps{{PSL(2,\mathbb C)}}
  \def\qed{{\hspace{2mm}{\small $\diamondsuit$}\goodbreak}}
  \def\pf{{\noindent{\bf Proof.\hspace{2mm}}}}
  \def\ni{{\noindent}}
  \def\sm{{{\mbox{\tiny M}}}}
   \def\sf{{{\mbox{\tiny F}}}}
   \def\sc{{{\mbox{\tiny C}}}}
  \def\ke{{\mbox{ker}(H_1(\partial M;\2)\ra H_1(M;\2))}}
  \def\et{{\mbox{\hspace{1.5mm}}}}

\title{Characterizing slopes for torus knots, II}

\author{{Yi NI}\\{\normalsize Department of Mathematics, Caltech}\\
{\normalsize 1200 E California Blvd, Pasadena, CA 91125}
\\{\small\it Email\/:\quad\rm yini@caltech.edu}
\\\\
{Xingru ZHANG}
\\
{\normalsize Department of Mathematics,
University at Buffalo}\\
{\small\it Email\/:\quad\rm xinzhang@buffalo.edu}}

\date{}
\maketitle

\begin{abstract}
A slope $\frac pq$ is called a characterizing slope for a given knot $K_0\subset S^3$ if whenever the $\frac pq$--surgery on a knot $K\subset S^3$ is homeomorphic to the $\frac pq$--surgery on $K_0$ via an orientation preserving homeomorphism, then $K=K_0$.
In a previous paper, we showed that, outside a certain finite set of slopes, only the negative integers could possibly be
 non-characterizing slopes for the torus knot $T_{5,2}$. Applying
    recent work of Baldwin--Hu--Sivek, we improve our result
    by showing that a nontrivial slope $\frac pq$ is a characterizing slope for $T_{5,2}$ if $\frac pq>-1$ and  $\frac pq\notin \{0,1, \pm\frac12,\pm\frac13\}$.
In particular, every nontrivial L-space slope of $T_{5,2}$ is characterizing for $T_{5,2}$. As a consequence, if a nontrivial  $\frac pq$-surgery on a non-torus knot in $S^3$ yields a manifold of finite fundamental group, then $|p|>9$.
\end{abstract}

\section{Introduction}

For  a knot $K\subset S^3$ and a slope  $\frac pq\in\mathbb Q\cup\{\infty\}$, let $S^3_{p/q}(K)$ be the manifold obtained by the $\frac pq$--surgery on $K$.
A slope $\frac pq$ is said to be {\it characterizing} for a given knot $K_0\subset S^3$ if whenever $S^3_{p/q}(K)\cong S^3_{p/q}(K_0)$ for a knot $K\subset S^3$, then $K=K_0$. Here, ``$\cong$'' stands for an orientation preserving homeomorphism. Obviously the trivial slope $\frac 10$ is never characterizing for any knot.

A long standing conjecture due to Gordon states that all nontrivial slopes are characterizing  for the unknot. This conjecture was originally proved using Monopole Floer homology \cite{KMOSz}, and there were also proofs via Heegaard Floer homology \cite{OSzGenus,OSzRatSurg}. Based on work of Ghiggini \cite{Gh}, Ozsv\'ath and Szab\'o \cite{OSz3141} proved the same result for the trefoil knot and the figure--$8$ knot. To date, the unknot, the trefoil knot and the figure--$8$ knot are the only knots for which it is known that all but finitely many slopes are characterizing.

The next simplest knot is the torus knot $T_{5,2}$. It is reasonable to expect that all nontrivial  slopes are characterizing for $T_{5,2}$. In \cite{NiZhang}, it is proved that a nontrivial slope  $\frac pq$ is characterizing for  $T_{5,2}$ unless $\frac pq$ is a negative integer or $-47\le p\le32$ and $1\le q\le8$. It is also known that the slopes \[8,9,10,11,12,\frac{17}2,\frac{19}2,\frac{21}2,\frac{23}2,\frac{28}3,\frac{29}3,\frac{31}3,\frac{32}3\] are characterizing \cite{BZ,Baker,GreeneCabling,NiZhangFinite}. Characterizing slopes for general torus knots have been studied in \cite{NiZhang,McCoySharp,McCoyTorus,NiCharSlope}.

In this paper, we further narrow down the range of possible non-characterizing slopes for $T_{5,2}$. Our main result is  the
following theorem.

\begin{thm}\label{thm:T52}
A nontrivial slope $\frac pq$ is  characterizing for $T_{5,2}$
if $\frac pq>-1$ and $\frac pq\notin\{0,1, \pm\frac12,\pm\frac13\}$.
\end{thm}

Recall that a rational homology $3$-sphere $Y$ is an {\it L-space} if the rank of $\widehat{HF}(Y)$ is equal to the order of $H_1(Y;\mathbb Z)$. For $T_{5,2}$, the surgery slopes which yield L-spaces are the ones greater than or equal to $3=2g(T_{5,2})-1$. Thus Theorem \ref{thm:T52} implies that all nontrivial L-space slopes of $T_{5,2}$ are characterizing for $T_{5,2}$.

\begin{cor}\label{cyclic-finite}If a nontrivial  $p/q$-surgery, with $|p|\leq 9$,  on a nontrivial knot $K$ in $S^3$ produces a manifold
of finite fundamental group, then $K$  is one of the torus knots
$T_{3,\pm2}$ and $T_{5,\pm2}$.
\end{cor}
\begin{proof}
By changing $K$ to its mirror image, we may also assume that $p/q$ is positive.

When $S^3_{p/q}(K)$ has cyclic fundamental group, i.e. when it is a lens space,
we may further assume that  $p/q=p$ is an integer for otherwise $K$ is a torus knot \cite{CGLS}
and $S^3_{p/q}(K)$ is never  a lens space when $p\leq 9$ and $q\geq 2$ \cite{Moser}.
Now  \cite[Theorem 1.4]{Greene}  implies that when $p\leq 9$,
the knot $K$ is a  fibered knot of genus at most $2$, and $S^3_p(K)\cong S^3_p(T_{3,2})$ or $S^3_p(T_{5,2})$.
Now we use \cite{Gh} or Theorem~\ref{thm:T52} to get that $K=T_{3,2}$ or $T_{5,2}$.

When $S^3_{p/q}(K)$ has finite but non-cyclic fundamental group,
it is shown in \cite[Theorems~2, 3 and Table~1]{Doig} that either
$K$ is $T_{3,2}$ or $T_{5,2}$, or $p/q=7$ or $8$,
$S^3_{p/q}(K)\cong S^3_{p/q}(T_{5,2})$.
We may now apply Theorem~\ref{thm:T52} to conclude that $K=T_{5,2}$.
\end{proof}

The bound  $9$  in Corollary~\ref{cyclic-finite} appears to be the best one could get for hyperbolic knots in $S^3$ with the current techniques. The $10$--surgery on $T_{4,3}$ is a spherical space form with non-cyclic fundamental group, and the current techniques could not rule out the possibility that the same surgery on a hyperbolic knot with the same knot Floer homology as $T_{4,3}$ yields the same manifold. Conjecturally on a hyperbolic knot $K$ in
$S^3$ if a nontrivial $\frac pq$-surgery yields a manifold of  finite fundamental group, then $|p|\geq 17$, a bound which  can be realized on the $(-2,3,7)$--pretzel knot.

Our proof of Theorem \ref{thm:T52}, given in the next  section, is mainly  based on
our earlier work \cite{NiZhang} and a  recent paper of Baldwin--Hu--Sivek \cite{BHS}.

\vspace{5pt}\noindent{\bf Acknowledgements.}\quad The first author was
partially supported by NSF grant
number DMS-1811900.


\section{Proof of Theorem \ref{thm:T52}}\label{sect:proof}

For a knot $K$ in $S^3$, $\D_K(T)$
denotes the symmetric Alexander polynomial of $K$.
The following theorem and remark are \cite[Theorem~4.1]{NiZhang} and the remark after the proof.

\begin{thm}\label{thm:Fibered}
Suppose  that $S^3_{p/q}(K)\cong S^3_{p/q}(T_{5,2})$ for a knot $K\subset S^3$ and a nontrivial slope $\frac pq$. Then one of the following two cases happens:

1) $K$ is a genus $(n+1)$ fibered knot for some $n\ge1$ with
\begin{equation}\label{eq:DeltaKn}
\Delta_K(T)=(T^{n+1}+T^{-(n+1)})-2(T^n+T^{-n})+(T^{n-1}+T^{1-n})+(T+T^{-1})-1.
\end{equation}

2) $K$ is a genus $1$ knot with $\Delta_K(T)=3T-5+3T^{-1}$.

Moreover,   if $$\frac pq\in\left\{\frac pq>1\right\}\cup \left\{
\frac pq<-6, |q|\ge2\right\},$$ then the number $n$ in the first case must be $1$.
\end{thm}

\begin{rem}\label{addendum}
We have the following addendum to Theorem~\ref{thm:Fibered}:
\newline
(a) If $p$ is even, then case 2) of Theorem~\ref{thm:Fibered} cannot happen and in case 1) of Theorem~\ref{thm:Fibered}, the number $n$ must be odd.
\newline
(b) If $p$ is divisible by $3$, then case 2) cannot happen and in case 1), the number $n$ is not divisible by $3$.
\end{rem}

The following result is implicitly contained in \cite[Subsection~4.1]{NiZhang}. Background information about Heegaard Floer homology can be found in \cite[Section~3]{NiZhang}.

\begin{prop}\label{prop:SameHFK}
Suppose that  $S^3_{p/q}(K)\cong S^3_{p/q}(T_{5,2})$ for a knot $K\subset S^3$ and a nontrivial slope  $\frac pq>1$. Then
\begin{equation}\label{eq:HFKisom}
\widehat{HFK}(S^3,K)\cong\widehat{HFK}(S^3, T_{5,2})
\end{equation} as a bigraded group.
\end{prop}
\begin{proof}
In \cite[Subsection~4.1]{NiZhang}, it is proved that $\Delta_K=\Delta_{T_{5,2}}$ and $V_0(K)=V_1(K)=1$. By \cite[Proposition~4.2]{NiZhang},
\[
t_s(K)=V_s(K)+\mathrm{rank}H_{\mathrm{red}}(A_s^+).
\]
Since $t_0(K)=t_1(K)=1$ and $t_s(K)=0$ when $s>1$, we have
\[
HF_{\mathrm{red}}(A_s^+)=0\quad\text{whenever }s\ge0.
\]
So $K$ is an L-space knot. Since $\Delta_K=\Delta_{T_{5,2}}$, we get (\ref{eq:HFKisom}) by \cite{OSzLens}.
\end{proof}

\begin{rem}
In \cite[Proposition~4.7]{NiZhang}, it is proved that if $\frac pq<-6$ and $|q|\ge2$, then $\Delta_K=\Delta_{T_{5,2}}$. However, we cannot conclude that (\ref{eq:HFKisom}) holds. Algebraically, it is possible that $\widehat{HFK}(S^3,K,1)$ has rank $3$, while $HF^+(S^3_{p/q}(K))\cong HF^+(S^3_{p/q}(T_{5,2}))$.
\end{rem}

For any hyperbolic knot $K$ in $S^3$,
a result of Gabai and Mosher \cite{Mosher} states that the
complement of $K$ contains an essential lamination $\lambda$ which has an associated degeneracy locus $d(\lambda)$  in the form of $d(\lambda)=\frac mn$, $(m,n)\in\mathbb Z^2\setminus\{(0,0)\}$, such that if \[\Delta(\frac pq, d(\lambda)):=|pn-qm|\ge2,\] then $\lambda$
remains an essential lamination in $S^3_{p/q}(K)$.
Since $S^3_{1/0}(K)=S^3$ does not contain any essential lamination, it follows that
$d(\lambda)=m/0$ or $m/1$.
Furthermore we have

{\bf Fact (i)} If $\Delta(\frac pq, d(\lambda))\ge3$, then $S^3_{p/q}(K)$ is an irreducible, atoroidal and non-Seifert fibered manifold \cite[Theorem~2.5]{Wu}.

{\bf Fact (ii)} If $K$ is fibered and $d(\lambda)=m/1$ where $\lambda$ is the stable lamination of $K$, then $|m|\ge2$  \cite[Theorem~8.8]{GabaiProblems}. (Another proof of this was given in \cite{Roberts2}.)

\begin{prop}\label{prop:p/q<1}
Suppose that  $S^3_{p/q}(K)\cong S^3_{p/q}(T_{5,2})$, $|\frac pq|<1$ and $\frac pq\notin\{0,\pm\frac12,\pm\frac13\}$. Then $K=T_{5,2}$.
\end{prop}

\begin{proof}
When $K$ is a torus knot, Proposition \ref{prop:p/q<1} holds by \cite[Lemma 4.2]{McCoyTorus}.

Next suppose that $K$ is a satellite knot.
Let $K'$ be a companion knot of $K$. Since $S^3_{p/q}(K)\cong S^3_{p/q}(T_{5,2})$ is atoroidal and irreducible (since $p/q\ne 10$) \cite{Moser}, the work of Gabai \cite{GSurg} implies that $K$ is a $0$--bridge or $1$--bridge braid in a tubular neighborhood $V$ of $K'$.
If $K$ is a $0$--bridge braid, then $K$ is a $(r, s)$-cable of $K'$
and the surgery slope $p/q$ must be $(qrs\pm 1)/q$ \cite{Go}.
If $K$ is a $1$--bridge braid, then it follows from \cite[Lemma 3.2]{G1bridge} that $\frac pq\in\mathbb Z$. In both cases we get a contradiction with the assumption that $|p/q|<1$.

So  $K$ is hyperbolic.
By Theorem~\ref{thm:Fibered}, either $K$ is fibered or $g(K)=1$.
If $K$ is fibered and $d(\lambda)=\frac m0$, since $|q|\ge3$, $\Delta(d(\lambda),\frac pq)\ge3$. Hence $S^3_{p/q}(K)$
is not Seifert fibered by Fact (i) above,  which contradicts  the assumption $S^3_{p/q}(K)\cong S^3_{p/q}(T_{5,2})$.
So we may assume that when $K$ is fibered,  $d(\lambda)=m/1$ for the stable lamination $\lambda$ of $K$, and therefore $|m|\geq 2$ by Fact (ii) above.
Since $|\frac pq|<1$ and $|m|\ge2$, $\Delta(d(\lambda),\frac pq)\ge3$. Again by Fact (i) we get a contradiction with the assumption
$S^3_{p/q}(K)\cong S^3_{p/q}(T_{5,2})$.

So $K$ is hyperbolic and $g(K)=1$. By  \cite[Theorem 1.5]{BGZ} $0<|p|\leq 3$
since $S^3_{p/q}(K)\cong S^3_{p/q}(T_{5,2})$ is a Seifert fibered space.
 By Remark~\ref{addendum}, $|p|=1$.
Since $\frac pq\notin\{0,\pm\frac12,\pm\frac13\}$, we again have $\Delta(d(\lambda),\frac pq)\ge3$ (whether $d(\lambda)=\frac m0$ or $\frac m1$), which leads to  a contradiction as in the preceding  paragraph.
\end{proof}

\begin{prop}\label{prop:p/q>1}
Suppose that  $S^3_{p/q}(K)\cong S^3_{p/q}(T_{5,2})$ for a nontrivial slope $\frac pq$ with $\frac pq>1$. Then $K=T_{5,2}$.
\end{prop}

\begin{proof}To get a contradiction we assume that $K\ne T_{5,2}$.

By Proposition \ref{prop:SameHFK}, we have $\widehat{HFK}(S^3,K)\cong\widehat{HFK}(S^3, T_{5,2})$.
In \cite[Section 3]{BHS}, Baldwin--Hu--Sivek proved that if $\widehat{HFK}(S^3,K)\cong\widehat{HFK}(S^3, T_{5,2})$ and $K\ne T_{5,2}$, then   $K$ is a hyperbolic, doubly-periodic, genus two fibered  knot,
 the degeneracy locus of the stable lamination of $K$ is $4$, and moreover
there exists a pseudo-Anosov $5$-braid $\beta$ whose closure $B=\hat{\beta}$ is an unknot with braid axis $A$, such that $K$ is the lift of $A$ in the branched double cover $\Sigma(S^3,B)\cong S^3$.
Let $V$ be the exterior of $A$ in $S^3$ and $M_K$ the exterior of $K$ in $S^3$. Then $V$ is a solid torus
and $M_K$ is a double branched cover of $V$ with
$B$ as the branched set in $V$.
Let $\t$ be the corresponding covering involution on $M_K$
 and $U$ the branched set in $M_K$. Then $U$ is the fixed point set of $\t$ which is a knot disjoint from $\p M_K$, and $(M_K,U)/\t=(V, B)$.
 The restriction of $\t$  on $\p M_K$  is a free action--an order two
rotation along the longitude factor of $\p M_K$.

We also use $M_K(p/q)$ to denote the surgery manifold $S^3_{p/q}(K)$
and similarly $V(p/q)$ for $S^3_{p/q}(A)$.
Note that  the involution $\t$ on $M_K$  extends to an  involution
$\t_{p/q}$ on $M_K(p/q)$.
In fact if we let $N_{p/q}$ denote the filling solid torus
in forming $M_K(p/q)=M_K\cup N_{p/q}$ and let $K_{p/q}$ be the center circle of $N_{p/q}$,
then  the fixed point set of $\t_{p/q}$ is
\begin{equation}\label{fixed point set for p/q}
Fix(\t_{p/q})=\left\{\begin{array}{l}U,\;\;\;\;\;\;\;\;\;\;\;\;\;\;\mbox{if $p$ is odd}
\\U\cup K_{p/q}, \;\;\;\mbox{if $p$ is even }\end{array}\right.
\end{equation}
and
\begin{equation}\label{branched set for p/q}
M_K(p/q)/\t_{p/q}=\left\{\begin{array}{l}V(p/2q),\;\;\;\;\;\;\;\mbox{if $p$ is odd}
\\V((p/2)/q), \;\;\;\mbox{if $p$ is even.}\end{array}\right.
\end{equation}

Let  $Y$ be the exterior of $U$ in $M_K$ and $W$ the exterior of $B$ in $V$. Then
$Y$ is a free double cover of $W$.
Since $B$ is the closure of a pseudo-Anosov braid in $V$,
$W$ is hyperbolic.
Hence $Y$ is also hyperbolic.

Note that $M_K(p/q)\cong M_{T_{5,2}}(p/q)$ is a Seifert fibered space whose base orbifold
is $S^2(2, 5, d)$ with $d=|p-10q|$ if $p/q\ne 10$, and is a connected sum of
two lens spaces of orders $2$ and $5$ if $p/q=10$ \cite{Moser}.
Since $K$ is a doubly periodic knot, by \cite{HS} $M_K(p/q)$ is irreducible and thus $p/q\ne 10$, by \cite{WZ} $M_K(p/q)$ is not a lens space, and by \cite{MM} $M_K(p/q)$ is not a prism manifold.
Thus $M_K(p/q)$ is a Seifert fibered space
whose base orbifold is $S^2(2,5, d)$ with $d=|p-10q|>2$.
So there is a unique Seifert structure on $M_K(p/q)$.
 By Thurston's Orbifold Theorem \cite{BP,CHK},
  we may assume that  the unique Seifert structure on $M_K(p/q)$
 is $\t_{p/q}$-invariant, i.e. $\t_{p/q}$ sends every Seifert fiber  to a Seifert fiber preserving the order of singularity.

Since the base orbifold of the Seifert fibered space $M_K(p/q)$
is orientable, the Seifert fibers of $M_K(p/q)$ can be  coherently oriented.
If $\t_{p/q}$ preserves the orientations of the Seifert fibers of $M_K(p/q)$, then
$Fix(\t_{p/q})$ consists of Seifert fibers (see \cite[Lemma 4.3]{BGZ}).
Since $K$ is hyperbolic, $K_{p/q}$ cannot be a component of $Fix(\t_{p/q})$.
Thus by Formula (\ref{fixed point set for p/q}),  $p$ is odd and $Fix(\t_{p/q})=U$, and by Formula (\ref{branched set for p/q}),
$M_K(p/q)/\t_{p/q}=V(p/2q)$.
Moreover, if we let  $Y(\p M_K, p/q)$ denote the Dehn filling of $Y$ along the  component
$\p M_K$ of $\p Y$ with the slope $p/q$, and similarly  define $W(\p V, p/2q)$,
then $Y(\p M_K, p/q)$ is Seifert fibered
and is a free double cover of
$W(\p V, p/2q)$. So the latter manifold $W(\p V, p/2q)$ is also Seifert fibered.
But $W(\p V, 1/0)$ is a solid torus (it is the exterior of the unknot $B$ in $S^3$).
We get a contradiction with \cite[Corollary 15]{MZ}.

Hence $\t_{p/q}$ reverses the orientations of the Seifert fibers of $M_K(p/q)$. Since $p/q>1$, we may assume that both $p$ and $q$ are positive.
Let $\mu$ be the meridian of $K$, then $[\mu]$ generates $H_1(M_K)\cong\mathbb Z$ and $H_1(M_K(p/q))\cong\mathbb Z/p\mathbb Z$. Clearly $\t_*[\mu]=[\mu]$, so
\begin{equation}\label{eq:InvId}
\t_*=\mathrm{id} \text{ on }H_1(M_K),\qquad (\t_{p/q})_*=\mathrm{id} \text{ on }H_1(M_K(p/q)).
\end{equation}

To apply a homological argument, we describe the Seifert fibered structure of $M_{T_{5,2}}(p/q)$ explicitly as follows.
Let $V_0\cup V_1$ be a standard  genus one Heegaard splitting of $S^3$. We may assume that  $T_{5,2}$ is embedded in $\partial V_0$
and is homologous to $5\mathcal L+2\mathcal M$, where $\mathcal L$ is the canonical longitude of $V_0$, and $\mathcal M$  the meridian of $V_0$. Let $\mu_0\subset \partial M_{T_{5,2}}$ be the meridian of $T_{5,2}$, let $\lambda_0$ be the canonical longitude of $T_{5,2}$ and
let $\mathcal C_i$ be the core of $V_i$, $i=1,2$. Then
 $M_{T_{5,2}}$ is Seifert fibered with $C_0$ and $C_1$ as two singular fibers of order $5$ and $2$ respectively.
 A regular fiber $\cal F$ of $M_{T_{5,2}}$ in $\p M_{T_{5,2}}$ has slope  $10\m_0+\l_0$.
If $p/q\ne 10$, the Seifert structure of $M_{T_{5,2}}$ extends to one on $M_{T_{5,2}}(p/q)$ such that
  the core $\mathcal C'$ of the filling solid torus is an order $d=|p-10q|$ singular fiber if $d>1$.
  In $H_1(M_{T_{5,2}}(p/q))$, we have
\begin{equation}\label{eq:SingFiber}
[\mathcal F]=10[\mu_0],\quad[\mathcal C_0]=2[\mu_0],\quad [\mathcal C_1]=5[\mu_0],\quad [\mathcal C']=\pm q'[\mu_0],
\end{equation}
where $q'\in\mathbb Z$ satisfies that $qq'\equiv1\pmod p$.

 Since $\t_{p/q}$ sends a regular fiber to a regular fiber reversing
its orientation, we see, using (\ref{eq:InvId}) and (\ref{eq:SingFiber}),  that $10[\mu_0]=-10[\mu_0]$ in $\mathbb Z/p\mathbb Z$. So $p|20$.
We claim that $p=1$ or $2$. Suppose otherwise that $p>2$.
We already know that $M_K(p/q)$ has three singular fibers of orders
$2, 5, d=|p-10q|>2$. If $d\ne 5$, then $\t_{p/q}$ sends the order $d$ singular fiber $\mathcal C'$ to itself with opposite orientation.
Hence $q'[\mu_0]=-q'[\mu_0]$ in $\mathbb Z/p\mathbb Z$. Since $\gcd(p,q')=1$, we get $p|2$, which is not possible. So we must have $d=5$, which, together with the condition $p|20$,  implies that $p/q=5$. By (\ref{eq:SingFiber}) the two order $5$ singular fibers $\mathcal C_0,\mathcal C'$ are homologous to $2[\mu_0]$ and $\pm[\mu_0]$ respectively. By (\ref{eq:InvId}), $\t_{p/q}$ must send $\mathcal C_0$ to itself, and $\mathcal C'$ to itself. So $2[\mu_0]=-2[\mu_0]$ in $\mathbb Z/5\mathbb Z$, which is not possible.

Recall that  $K$ has a degeneracy locus $4$. Since $\frac pq\le2$, we have $\Delta(\frac pq, 4)\ge3$ unless $\frac pq=2$.
By Fact (i), we only need to consider the case $\frac pq=2$.

By Formula (\ref{fixed point set for p/q}) the fixed point set of
$\t_2$ is $U\cup K_2$ and by  Formula (\ref{branched set for p/q}),
 $M_K(2)/\t_{2}=V(1)$ which is $S^3$. Hence the branched set $B\cup K_{2}^*$ in
 $M_K(2)/\t_{2}=V(1)=S^3$
is a Montesinos link of two components \cite{Mon}, where $K_2^*$ is the image of $K_2$ under the map $M_K(2)\ra V(1)$, which is also the core of the filling solid torus of $V(1)$.
Note  that $K_{2}^*$ is an unknot in $S^3$
while  $B$ is the closure of a $5$-braid in the exterior of $K_{2}^*$ which is a solid torus.
Hence the linking number between $B$ and $K^*_{2}$ is $5$.

\begin{figure}[!ht]
\centerline{\includegraphics{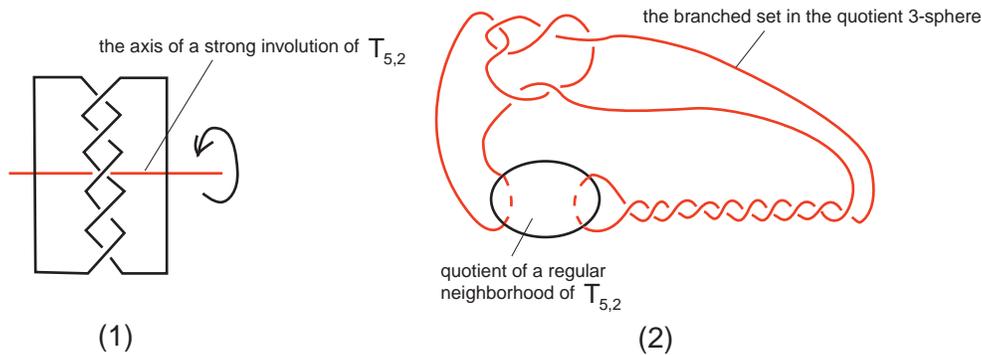}} \caption{The strong involution on $(S^3, T_{5,2})$
and the quotient}\label{involution on (5,2)}
\end{figure}

\begin{figure}[!ht]
\centerline{\includegraphics{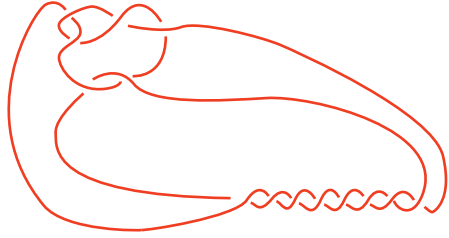}} \caption{The branched sets for the $2$-surgery on $T_{5,2}$}\label{branched set for slope 2}
\end{figure}

On the other hand  $M_{T_{5,2}}(2)$ is a double branched cover of $S^3$ whose branched set in $S^3$ can be explicitly constructed
 as follows.
 The knot $T_{5,2}$ is strongly invertible
and the quotient spaces under the strong involution of $S^3$, a regular neighborhood of $T_{5,2}$
 and the exterior, as well as the branched  set
are as shown in Figure~\ref{involution on (5,2)}.
Furthermore  $M_{T_{5,2}}(2)$ is a double branched cover of
$S^3$ whose branch set can be obtained by replacing the rational  $1/0$-tangle in Figure~\ref{involution on (5,2)}~(2)
with the rational $2$-tangle.
The resulting branched set is the  two-component
Montesinos link shown in Figure~\ref{branched set for slope 2}.

Note that   $M_K(2)=M_{T_{5,2}}(2)$ has base orbifold $S^2(2,5,8)$ and  so it has a unique Seifert fibration structure.
Since the number of singular fibers is $3$, the classification of Montesinos links \cite{BurdeZieschang} implies that $M_{T_{5,2}}(2)$ is the double branched cover of $S^3$
over a unique  Montesinos link with $3$ rational tangles.
Therefore the link 
 shown in Figure~\ref{branched set for slope 2} 
 should be the link $B\cup K_2^*$.
 However the linking number between the two components of the link in Figure~\ref{branched set for slope 2} is $3$, yielding a final contradiction with the early conclusion that the linking number between $B$ and $K_2^*$ is $5$.
\end{proof}

\begin{rem}
Proposition \ref{prop:p/q>1} can be proved without using the
degeneracy locus condition.
\end{rem}

Now the combination of Propositions \ref{prop:p/q<1} and \ref{prop:p/q>1} gives Theorem \ref{thm:T52}.


\begin{thebibliography}{H}

\bib{Baker}{article}{
   author={Baker, Kenneth L.},
   title={Small genus knots in lens spaces have small bridge number},
   journal={Algebr. Geom. Topol.},
   volume={6},
   date={2006},
   pages={1519--1621},
   issn={1472-2747},
}

\bib{BHS}{article}{
 author={Baldwin, John A.},
 author={Hu, Ying},
   author={Sivek, Steven},
   title={Khovanov homology and the cinquefoil},
   status={preprint},
   date={2021},
   eprint={https://arxiv.org/abs/2105.12102}
}


\bib{BP}{article}{
   author={Boileau, Michel},
   author={Porti, Joan},
   title={Geometrization of 3-orbifolds of cyclic type},
   language={English, with English and French summaries},
   note={Appendix A by Michael Heusener and Porti},
   journal={Ast\'{e}risque},
   number={272},
   date={2001},
   pages={208},
}

\bib{BGZ}{article}{
   author={Boyer, Steven},
   author={Gordon, Cameron McA.},
   author={Zhang, Xingru},
   title={Dehn fillings of knot manifolds containing essential
   once-punctured tori},
   journal={Trans. Amer. Math. Soc.},
   volume={366},
   date={2014},
   number={1},
   pages={341--393},
   issn={0002-9947},
}

\bib{BZ}{article}{
   author={Boyer, S.},
   author={Zhang, X.},
   title={Finite Dehn surgery on knots},
   journal={J. Amer. Math. Soc.},
   volume={9},
   date={1996},
   number={4},
   pages={1005--1050},
   issn={0894-0347},
}

\bib{BurdeZieschang}{book}{
   author={Burde, Gerhard},
   author={Zieschang, Heiner},
   title={Knots},
   series={De Gruyter Studies in Mathematics},
   volume={5},
   edition={2},
   publisher={Walter de Gruyter \& Co., Berlin},
   date={2003},
   pages={xii+559},
}	


\bib{CHK}{book}{
   author={Cooper, Daryl},
   author={Hodgson, Craig D.},
   author={Kerckhoff, Steven P.},
   title={Three-dimensional orbifolds and cone-manifolds},
   series={MSJ Memoirs},
   volume={5},
   note={With a postface by Sadayoshi Kojima},
   publisher={Mathematical Society of Japan, Tokyo},
   date={2000},
   pages={x+170},
}

\bib{CGLS}{article}{
   author={Culler, Marc},
   author={Gordon, C. McA.},
   author={Luecke, J.},
   author={Shalen, Peter B.},
   title={Dehn surgery on knots},
   journal={Ann. of Math. (2)},
   volume={125},
   date={1987},
   number={2},
   pages={237--300},
}

\bib{Doig}{article}{
   author={Doig, Margaret I.},
   title={Finite knot surgeries and Heegaard Floer homology},
   journal={Algebr. Geom. Topol.},
   volume={15},
   date={2015},
   number={2},
   pages={667--690},
}


\bib{GSurg}{article}{
   author={Gabai, David},
   title={Surgery on knots in solid tori},
   journal={Topology},
   volume={28},
   date={1989},
   number={1},
   pages={1--6},
   issn={0040-9383},
}

\bib{G1bridge}{article}{
   author={Gabai, David},
   title={$1$-bridge braids in solid tori},
   journal={Topology Appl.},
   volume={37},
   date={1990},
   number={3},
   pages={221--235},
}

\bib{GabaiProblems}{article}{
   author={Gabai, David},
   title={Problems in foliations and laminations},
   conference={
      title={Geometric topology},
      address={Athens, GA},
      date={1993},
   },
   book={
      series={AMS/IP Stud. Adv. Math.},
      volume={2},
      publisher={Amer. Math. Soc., Providence, RI},
   },
   date={1997},
   pages={1--33},
}



\bib{Gh}{article}{
   author={Ghiggini, Paolo},
   title={Knot Floer homology detects genus-one fibred knots},
   journal={Amer. J. Math.},
   volume={130},
   date={2008},
   number={5},
   pages={1151--1169},
}

\bib{Go}{article}{
   author={Gordon, C. McA.},
   title={Dehn surgery and satellite knots},
   journal={Trans. Amer. Math. Soc.},
   volume={275},
   date={1983},
   number={2},
   pages={687--708},
}



\bib{Greene}{article}{
   author={Greene, Joshua Evan},
   title={The lens space realization problem},
   journal={Ann. of Math. (2)},
   volume={177},
   date={2013},
   number={2},
   pages={449--511},
   issn={0003-486X},
}


\bib{GreeneCabling}{article}{
   author={Greene, Joshua Evan},
   title={L-space surgeries, genus bounds, and the cabling conjecture},
   journal={J. Differential Geom.},
   volume={100},
   date={2015},
   number={3},
   pages={491--506},
   issn={0022-040X},
}


\bib{HS}{article}{
   author={Hayashi, Chuichiro},
   author={Shimokawa, Koya},
   title={Symmetric knots satisfy the cabling conjecture},
   journal={Math. Proc. Cambridge Philos. Soc.},
   volume={123},
   date={1998},
   number={3},
   pages={501--529},
}

\bib{KMOSz}{article}{
   author={Kronheimer, P.},
   author={Mrowka, T.},
   author={Ozsv\'{a}th, P.},
   author={Szab\'{o}, Z.},
   title={Monopoles and lens space surgeries},
   journal={Ann. of Math. (2)},
   volume={165},
   date={2007},
   number={2},
   pages={457--546},
}



\bib{McCoyTorus}{article}{
   author={McCoy, Duncan},
   title={Non-integer characterizing slopes for torus knots},
   journal={Comm. Anal. Geom.},
   volume={28},
   date={2020},
   number={7},
   pages={1647--1682},
   issn={1019-8385},
}


\bib{McCoySharp}{article}{
   author={McCoy, Duncan},
   title={Surgeries, sharp 4-manifolds and the Alexander polynomial},
   status={preprint},
   date={2014},
   eprint={https://arxiv.org/abs/1412.0572}
}

\bib{MZ}{article}{
   author={Menasco, W.},
   author={Zhang, X.},
   title={Notes on tangles, 2-handle additions and exceptional Dehn
   fillings},
   journal={Pacific J. Math.},
   volume={198},
   date={2001},
   number={1},
   pages={149--174},
}

\bib{MM}{article}{
   author={Miyazaki, Katura},
   author={Motegi, Kimihiko},
   title={Seifert fibering surgery on periodic knots},
   booktitle={Proceedings of the First Joint Japan-Mexico Meeting in
   Topology (Morelia, 1999)},
   journal={Topology Appl.},
   volume={121},
   date={2002},
   number={1-2},
   pages={275--285},
}


\bib{Mon}{article}{
   author={Montesinos, Jos\'{e} M.},
   title={Seifert manifolds that are ramified two-sheeted cyclic coverings},
   language={Spanish},
   journal={Bol. Soc. Mat. Mexicana (2)},
   volume={18},
   date={1973},
   pages={1--32},
}



\bib{Moser}{article}{
   author={Moser, Louise},
   title={Elementary surgery along a torus knot},
   journal={Pacific J. Math.},
   volume={38},
   date={1971},
   pages={737--745},
}

\bib{Mosher}{article}{
   author={Mosher, Lee},
   title={Laminations and flows transverse to finite depth foliations. Part I:  branched surfaces and dynamics},
   status={preprint},
}


\bib{NiFibred}{article}{
   author={Ni, Yi},
   title={Knot Floer homology detects fibred knots},
   journal={Invent. Math.},
   volume={170},
   date={2007},
   number={3},
   pages={577--608},
}

\bib{NiCharSlope}{article}{
   author={Ni, Yi},
   title={Seifert fibered and reducible surgeries on hyperbolic fibered knots},
   status={preprint},
   date={2020},
   eprint={https://arxiv.org/abs/2007.11774}
}

\bib{NiZhang}{article}{
   author={Ni, Yi},
   author={Zhang, Xingru},
   title={Characterizing slopes for torus knots},
   journal={Algebr. Geom. Topol.},
   volume={14},
   date={2014},
   number={3},
   pages={1249--1274},
}

\bib{NiZhangFinite}{article}{
   author={Ni, Yi},
   author={Zhang, Xingru},
   title={Finite Dehn surgeries on knots in $S^3$},
   journal={Algebr. Geom. Topol.},
   volume={18},
   date={2018},
   number={1},
   pages={441--492},
   issn={1472-2747},
}



\bib{OSzGenus}{article}{
   author={Ozsv\'{a}th, Peter},
   author={Szab\'{o}, Zolt\'{a}n},
   title={Holomorphic disks and genus bounds},
   journal={Geom. Topol.},
   volume={8},
   date={2004},
   pages={311--334},
   issn={1465-3060},
}

\bib{OSzLens}{article}{
   author={Ozsv\'{a}th, Peter},
   author={Szab\'{o}, Zolt\'{a}n},
   title={On knot Floer homology and lens space surgeries},
   journal={Topology},
   volume={44},
   date={2005},
   number={6},
   pages={1281--1300}
}



\bib{OSzRatSurg}{article}{
   author={Ozsv\'{a}th, Peter},
   author={Szab\'{o}, Zolt\'{a}n},
   title={Knot Floer homology and rational surgeries},
   journal={Algebr. Geom. Topol.},
   volume={11},
   date={2011},
   number={1},
   pages={1--68},
   issn={1472-2747},
}

\bib{OSz3141}{article}{
   author={Ozsv\'{a}th, Peter},
   author={Szab\'{o}, Zolt\'{a}n},
   title={The Dehn surgery characterization of the trefoil and the figure
   eight knot},
   journal={J. Symplectic Geom.},
   volume={17},
   date={2019},
   number={1},
   pages={251--265},
}



\bib{Roberts2}{article}{
   author={Roberts, Rachel},
   title={Taut foliations in punctured surface bundles. II},
   journal={Proc. London Math. Soc. (3)},
   volume={83},
   date={2001},
   number={2},
   pages={443--471},
   issn={0024-6115},
}



\bib{WZ}{article}{
   author={Wang, Shi Cheng},
   author={Zhou, Qing},
   title={Symmetry of knots and cyclic surgery},
   journal={Trans. Amer. Math. Soc.},
   volume={330},
   date={1992},
   number={2},
   pages={665--676},
}

\bib{Wu}{article}{
   author={Wu, Ying-Qing},
   title={Sutured manifold hierarchies, essential laminations, and Dehn
   surgery},
   journal={J. Differential Geom.},
   volume={48},
   date={1998},
   number={3},
   pages={407--437},
}
		



\end{thebibliography}
\end{document}